\documentstyle[epsf,12pt]{article}
\normalbaselines
\catcode `\@=11
\@addtoreset{equation}{section}

\def\section{\@startsection {section}{1}{\z@}{-3.5ex plus -1ex minus
     -.2ex}{2.3ex plus .2ex}{\normalsize\bf}}
\def\subsection{\@startsection{subsection}{2}{\z@}{-3.25ex plus -1ex minus
 -.2ex}{1.5ex plus .2ex}{\normalsize\bf}}

\def\thebibliography#1{\section*{References\markboth
  {REFERENCES}{REFERENCES}}\list
  {[\arabic{enumi}]}{\settowidth\labelwidth{[#1]}\leftmargin\labelwidth
  \advance\leftmargin\labelsep
  \usecounter{enumi}}
  \def\newblock{\hskip .11em plus .33em minus -.07em}
  \sloppy
  \sfcode`\.=1000\relax}
 
\catcode `\@=12

\begin{document}

\begin{flushright}
Preprint DM/IST 9/97
\end{flushright}
\vspace*{1cm}
\begin{center}
{ \bf REFLECTIONS ON TOPOLOGICAL QUANTUM FIELD THEORY
\footnote{Talk presented at the XVth Workshop on Geometric Methods in 
Physics, 
{\it Quantizations, Deformations and Coherent States}, 
Bia{\l}owie\.{z}a, Poland, July 1 - 7, 1996. To appear in {\it Rep.  Math.
Phys.}} 
}\vspace{1cm}\\ \end{center}

\begin{center}
Roger Picken$^{1}$\\ \vspace{0.3cm}
${}^{1}$  Departamento de Matem\'{a}tica \& Centro de Matem\'{a}tica 
Aplicada,\\
Instituto Superior T\'{e}cnico, Av. Rovisco Pais,\\
1096 Lisboa Codex, Portugal.
\end{center}

\vspace*{0.5cm}

\begin{abstract}
\noindent
The aim of this article is to introduce some basic notions of Topological 
Quantum
Field Theory (TQFT) and to consider a modification of TQFT, applicable to
embedded manifolds. After an introduction based around a simple example
(Section 1) the notion of a $d$-dimensional TQFT is defined in
category-theoretical terms, as a certain type of functor from a category of
$d$-dimensional cobordisms to the category of vector spaces (Section 2). A
construction due to Turaev,  an operator-valued invariant of tangles, is
discussed in Section 3. It bears a strong resemblance to 1-dimensional 
TQFTs,
but carries much richer structure due to the fact that the 1-dimensional
manifolds involved are embedded in a 3-dimensional space. This leads us, in
Section 4, to propose a class of TQFT-like theories, appropriate to 
embedded,
rather than pure, manifolds.
\end{abstract}

\section{\hspace{-4mm}.\hspace{2mm}INTRODUCTION}

\hspace*{0.8cm}
To introduce the idea of TQFT, we start by considering the simplest compact
phase space, namely $S^2$ with symplectic form $\omega$ given in terms of a
local complex coordinate $z$ by
$
\omega = idz\wedge d\bar{z}/(1+|z|^2)^2
$.
Locally $\omega$ is exact, $\omega = d\theta$, where $\theta$ is a 1-form, 
but
globally $\omega$ is the curvature of a non-trivial complex line bundle 
(i.e.
having fibre isomorphic to $\bf C$) with connection $\theta$. In geometric
quantization, the Hilbert space corresponding to $(S^2, \omega)$ is the 
space
of holomorphic sections of this line bundle. It can easily be shown that 
there
are only two linearly independent holomorphic sections, and thus the Hilbert
space is isomorphic to ${\bf C}^2$. (This illustrates the general principle
that a finite phase space volume gives rise to a finite Hilbert space
dimension.)

Turning to dynamics, since we are describing a topological theory the
Hamiltonian operator $\hat{H}$ is taken to be simply the zero $2\times 2$ 
matrix. 
Thus the quantum evolution operator $\exp (i\hat{H}t)$ is the
$2\times 2$ identity matrix and its trace is $2$.
We will now associate 1-dimensional manifolds and operators in the following
way: the interval I is associated with $\exp (i\hat{H}t)$ and the circle 
$S^1$
with ${\rm Tr}\, \exp (i\hat{H}t)$ (regarded as an operator from $\bf C$ 
to $\bf C$).
The idea is that the path integral representation of each operator is 
given by
an integral over fields defined on the corresponding $(0+1)$-dimensional
spacetime manifold. Incidentally, from the Legendre transformation one finds
that the action is $S=\int \theta$, since $0=H=\theta - {\cal L}$ where 
$\cal
L$ is the Lagrange density. Thus the factor $\exp iS$ in the path integral
 may
be interpreted geometrically as the holonomy of the connection 
$\theta$ \cite{at}.
\newpage 
Irrespective of how the above assignment is motivated, from a mathematical
standpoint it distinguishes between the two fundamental 1-dimensional 
manifolds
I and $S^1$ and thus provides a topological classification of (connected)
1-manifolds. Of course in one dimension there is very little to classify, 
but
the same feature persists for TQFTs relating to higher-dimensional manifolds, where the
classification problem can pose very great challenges.

In the next section we shall describe an elegant axiomatic approach to TQFT,
which allows one to bypass the difficulties associated with path 
integrals and captures
certain essential features common to all TQFTs.

\section{\hspace{-4mm}.\hspace{2mm}THE DEFINITION OF TQFT}

\hspace*{0.8cm}
One of the main approaches to defining the notion of a TQFT uses 
the elegant language of
category theory. For this reason we start with a very brief 
discussion of categories and
functors.

A {\em category} $\cal C$ consists of
(a) a class of {\em objects},
(b) for any ordered pair of objects $(V,W)$ a set ${\rm Hom}(V,W)$ of
{\em morphisms} $V\stackrel{f}{\rightarrow}W$, and 
(c) an operation of {\em composition} of morphisms,
$(V\stackrel{f}{\rightarrow}W,\, U\stackrel{g}{\rightarrow}V )
\mapsto\,\,U\stackrel{f\circ g}{\longrightarrow}W$
satisfying,
(C1) {\em associativity} of composition, and
(C2) existence of an {\em identity morphism} for each object (i.e. for
any $X$ there exists a morphism $X\stackrel{{\bf 1}_X}{\longrightarrow}X$ 
such that,
for any 
$X\stackrel{f}{\longrightarrow}Y$, 
$f\circ
{\bf 1}_X = f$ and ${\bf 1}_Y\circ f=f$).

Intuitively, the objects are sets or spaces, possibly endowed with 
additional
structures, and the morphisms are structure-preserving maps. Two examples
are:

\vglue 2mm
\noindent 1. ${\bf Vect}(k)$, the category whose objects are all 
finite-dimensional vector spaces over the ground ring $k$, and 
whose morphisms are
all $k$-linear maps.
\vglue 2mm
\noindent 2. $d-{\bf Cobord}$, the category whose objects are smooth, 
compact,
oriented $(d-1)$-dimensional manifolds without boundary, and whose 
morphisms are
{\em cobordisms} i.e. smooth, compact, oriented $d$-dimensional manifolds 
with boundary,
identified up to orientation-preserving diffeomorphisms which restrict 
to the
identity on the boundary. 
\vglue 2mm

\begin{figure}[h]
\hspace*{3cm}
\epsfysize=2.4in
\epsfbox{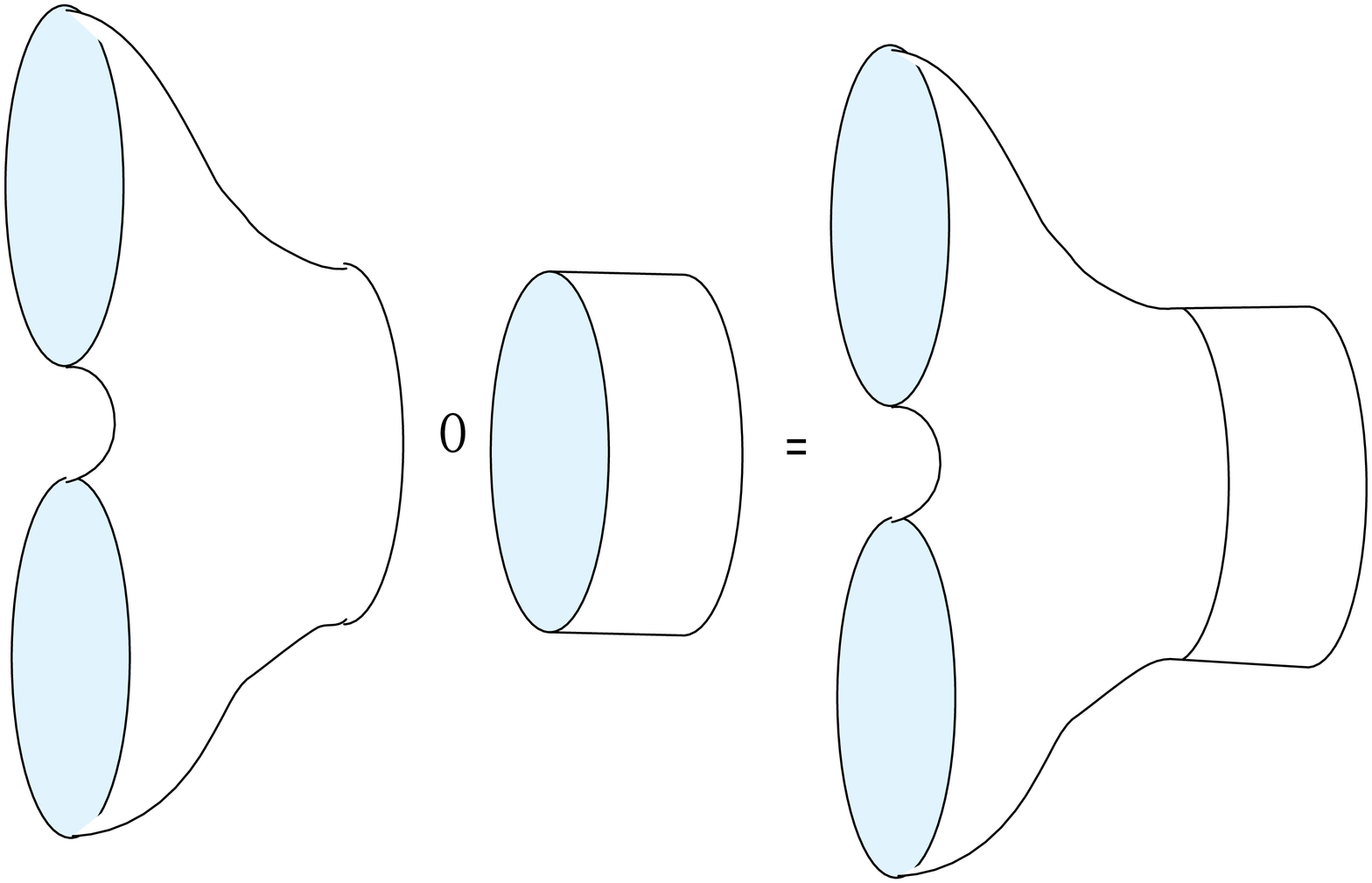}
\caption{Composition of 2-cobordisms}
\label{fig:trousers}
\end{figure}

The unusual feature of this second example is that the morphisms are 
not maps between
manifolds in the usual sense, but are themselves (equivalence classes of)
 manifolds of
one dimension higher than the the objects. The middle surface in
Figure~\ref{fig:trousers} depicts a 2-cobordism from the circle $S^1$ 
to the disjoint
union of two copies of the circle, namely the familiar ``trousers'' surface.

For simplicity, in what follows we will
describe cobordisms as
$d$-dimensional manifolds, rather than equivalence classes thereof.  For
$\Sigma_1$ and
$\Sigma_2$ two objects, a cobordism from $\Sigma_1$ to $\Sigma_2$, written
$\Sigma_1\stackrel{M}{\longrightarrow}
\Sigma_2$, is a manifold $M$ whose boundary is $\overline{\Sigma}_1\amalg
\Sigma_2$, where $\overline{\Sigma}_1$ means $\Sigma_1$ with the opposite
orientation and $\amalg$ is the disjoint union. The composition of a 
cobordism from 
$\Sigma_1$ to $\Sigma_2$ and one from $\Sigma_2$ to $\Sigma_3$ is given 
by gluing the
two cobordisms along $\Sigma_2$. 
The identity morphism for an
object $\Sigma$ is the cylinder cobordism $\Sigma\times I$, where $I$ 
is the
standard interval. Figure~\ref{fig:trousers} depicts the composition of 
the identity
morphism for the circle and the ``trousers'' cobordism. 

A  fundamental concept in category theory is that of a {\em functor} 
between
two categories: given two categories $\cal C$ and $\cal C'$, a (covariant) 
{\em
functor} from 
$\cal C$ to $\cal C'$ is
(a) an assigment to each object $V$ of $\cal C$ of an object $F(V)$ of
$\cal C'$, and 
(b) an assigment to each morphism $V\stackrel{f}{\longrightarrow}W$ of
$\cal C$ of a morphism
$F(V)\stackrel{F(f)}{\longrightarrow}F(W)$ of $\cal C'$,
such that
(F1) $F({\bf 1}_V)={\bf 1}_{F(V)}$ for every object $V$ of $\cal C$, and
(F2) $F(f\circ g)=F(f)\circ F(g)$ for every pair of morphisms $f$, $g$
of $\cal C$ whose composition is defined.

Categories may possess various kinds of additional structures, of which 
we mention
the following:
(A) a {\em product} of objects, and corresponding product of morphisms,
(B) a {\em unit object} and {\em unit morphism} with respect to this
product, and 
(C) an {\em involution} on the class of objects.
The categories ${\bf Vect}(k)$ and $d-{\bf Cobord}$ introduced above possess
all three structures.
In ${\bf Vect}(k)$
(A) is the tensor product $(V,W)\mapsto V\otimes W$  together with the 
tensor product of
two linear maps
$(V\stackrel{f}{\rightarrow}W, V'\stackrel{g}{\rightarrow}W')\mapsto 
(V\otimes V'\stackrel{f\otimes g}{\longrightarrow}W\otimes W')$,
(B) is the ring $k$ ($k\otimes V=V$,
etc.) together with
$1$ (regarded as a linear map
$k\rightarrow k$), and finally
(C) is $V\mapsto V^\ast$ (the dual vector space of $V$).
(Here we are identifying $V$ and $V^{\ast\ast}$.)
In $d-{\bf Cobord}$
(A)  is the disjoint union
$(\Sigma_1,\Sigma_2)\mapsto
\Sigma_1\amalg
\Sigma_2$ and 
$(\Sigma_1\stackrel{M}{\longrightarrow}{\Sigma_1}',
\Sigma_2\stackrel{N}{\longrightarrow}{\Sigma_2}')\mapsto  \Sigma_1\amalg
\Sigma_2\stackrel{M\amalg N}{\longrightarrow}{\Sigma_1}'\amalg {\Sigma_2}'$,
(B)  is $\emptyset$, the empty
$(d-1)$-dimensional manifold  
($\emptyset
\amalg
\Sigma = \Sigma$ etc.) together with $\emptyset$, the empty
$d$-dimensional manifold, and
(C) is $\Sigma \mapsto \overline{\Sigma}$ (the
manifold
$\Sigma$ with the opposite orientation).

We are now in a position to define what a TQFT is:
{\em A ($d$-dimensional) topological quantum field theory  is a functor
 $d-{\bf
Cobord}\longrightarrow {\bf Vect}(k)$ respecting the structures (A), (B) 
and (C).}

A TQFT is commonly described by assignments $\Sigma \mapsto V_\Sigma$ 
for objects
and $M\mapsto Z_M$ for morphisms. The definition implies that these 
obey the
 properties
(F1) $Z_{\Sigma\times I}={\rm id}_{V_\Sigma}$,
(F2) $Z_{M\circ M'}=Z_M\circ Z_{M'}$,
(AA) $V_{\Sigma\amalg{\Sigma}'}=V_{\Sigma}\otimes V_{{\Sigma}'}$ and
$Z_{M\amalg M'}=Z_M \otimes Z_{M'}$,
(BB) $V_\emptyset = k$; $Z_\emptyset = 1$, and
(CC) $V_{\overline{\Sigma}}={V_{\Sigma}}^\ast$.

As mentioned in the previous section, the physical intuition underlying 
this definition is
that
$V_\Sigma$ is the Hilbert space associated with the spacelike manifold 
$\Sigma$ and $Z_M$
is the path integral associated to the spacetime $M$.

It should be mentioned at this stage that there are other approaches to 
defining a
TQFT. In particular, the original axiomatic definition of TQFT due to 
Atiyah \cite{at},
has a slightly different flavour. In that approach a TQFT is an 
assignment  $\Sigma
\mapsto V_\Sigma$ of vector spaces to closed $(d-1)$-manifolds, 
together with an
assignment $M\mapsto Z_M$, which assigns to a $d$-dimensional 
manifold $M$, with
boundary $\partial M=\Sigma$, an {\em element} of the vector space 
$V_\Sigma$,
both assignments obeying a number of axioms. The definition in terms 
of a functor from
$d-{\bf Cobord}$ to ${\bf Vect}(k)$ has been used by various authors, 
e.g.
\cite{yetter}\cite{baezdolan}\cite{crane}\cite{sawin1}.

To give an idea of how this axiomatic definition of TQFT works, we will 
reanalyse
the TQFT discussed in the introduction from this standpoint. This is a 
1-dimensional
TQFT i.e. a functor from $1-{\bf Cobord}$ to ${\bf Vect}({\bf C})$. The
objects of
$1-{\bf Cobord}$, being 0-dimensional manifolds, are either $\emptyset$ 
or the
disjoint union of single points. Thus, by the properties (AA) and (BB) 
above, it is
enough to specify $V_p$, where $p$ stands for a point, and we set 
$V_p={\bf C}^2$.
 Turning to the morphisms, connected 1-manifolds are
diffeomorphic to either the interval $I$ or the circle $S^1$. By (F1) 
above $Z_I={\rm
id}_{{\bf C}^2}$. Thus to complete the description of this TQFT we need 
to obtain
$Z_{S^1}$. Now $\partial S^1=\emptyset$ so $S^1$ is a cobordism from 
$\emptyset$ to
$\emptyset$, and under the TQFT functor is transported to a linear map
$Z_{S^1}:V_\emptyset = {\bf C}\longrightarrow V_\emptyset = {\bf C}$.  
Thus $Z_{S^1}$ may
be identified with a number and it remains to show that this number is 
2, i.e. the trace
of the $2\times 2$ identity matrix.  
This is a relatively simple exercise, but for reasons of space we will 
not go into
details.

The discussion of TQFT in this section and the previous one were mainly 
intended to
convey a flavour of the subject. It is clear that the example of TQFT 
which was
studied has very little interest for the topological classification of 
manifolds,
since in dimension 1 the only (connected) manifolds to classify are the 
interval and
the circle. One way of achieving richer results is to go up in dimension. 
A notable
example is of course the 3-manifold invariant arising from a 3-dimensional 
TQFT based
on the Chern-Simons action \cite{qft+jones}. It would go too far to 
discuss here the many
other examples of higher dimensional TQFTs, since this is a vast and 
expanding
subject. Thus we limit ourselves to giving some references to the 
literature
\cite{baezdolan}\cite{crane}\cite{turbook}\cite{blauthomp}\cite{masbaum}
\cite{sawin2}\cite{freed}.
In the next two sections however we will propose a second way to achieve 
richer structure
for TQFTs, namely by considering a modified cobordism category whose 
objects and
morphisms are embedded manifolds. 

\section{\hspace{-4mm}.\hspace{2mm}OPERATOR INVARIANTS OF TANGLES}

\begin{figure}[h]
\hspace*{4cm}
\epsfysize=2in\epsffile{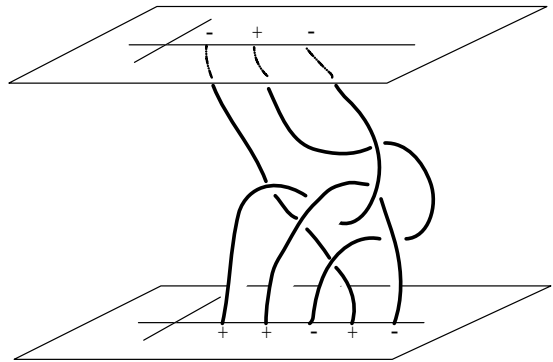}
\caption{A tangle}
\label{fig:tangle}
\end{figure}

\hspace*{0.8cm}
In this section we will look at a motivating example of the kind of 
``embedded''
TQFT alluded to at the end of the previous section. The construction 
of tangle
invariants to be described is due to Turaev \cite{tur}. We will present 
tangles as
a category, called ${\bf OTa}$ by Turaev (oriented tangles). The objects 
of this
category  are finite sequences of $+$ or $-$ signs, interpreted 
geometrically
as oriented points lying in integer positions on the positive $x_1$ 
axis in
${\bf R}^2$ coordinatised by $(x_1,x_2)$. The morphisms of ${\bf OTa}$ are
(oriented) tangles: consider ${\bf R}^3$ coordinatised by $(x_1,x_2,x_3)$ 
and
two parallel planes in ${\bf R}^3$, being $x_3=0$ and $x_3 =1$. 
Both planes are
to be thought of as objects in the above sense, and thus contain a 
number of
oriented points along the lines $x_2=x_3=0$ and $x_2=0,\, x_3=1$ 
respectively. A
tangle is a 1-dimensional oriented submanifold of
${\bf R}^2\times[0,1]\subset {\bf R}^3$ whose boundary consists of 
the oriented
points in the top and bottom planes, these being the only points 
where the
tangle intersects the top and bottom planes. Furthermore two 
tangles are
identified if one is carried to the other by an orientation-preserving
diffeomorphism of ${\bf R}^2\times [0,1]$ which is the identity 
restricted to
the top and bottom planes (see Figure~\ref{fig:tangle}).

A tangle is regarded as a morphism from the orientation-reversed object 
in the 
lower plane to the object in the upper plane, in the same way as a 
cobordism
from $\Sigma_1$ to $\Sigma_2$ was a manifold  $M$ with boundary 
$\overline{\Sigma}_1\amalg
\Sigma_2$. Thus the tangle in Figure~\ref{fig:tangle} is a morphism from
$(--+-+)$ to $(-+-)$. (Here we have adopted a different convention to
\cite{tur}, where points associated to upwards-pointing strands are 
labelled $+$ and
to downwards-pointing strands $-$. Our convention
aims to make contact with the TQFT axioms from the previous section.)

If $T_1$ is a tangle from $O_1$ to $O_2$ and $T_2$ is a tangle from $O_2$ to
$O_3$, we form the composition $T_2\circ T_1$ by concatenation, placing the
second tangle on top of the first and shrinking in the $x_3$ direction. The
identity morphism for an object is simply the corresponding trivial tangle,
which has all its strands parallel to the $x_3$ axis, oriented up or down as
appropriate.

The category ${\bf OTa}$ also possesses extra structures analogous to 
those of
$d-{\bf Cobord}$, namely (A) a product $\otimes$ of objects and
corresponding product of morphisms/tangles, given by juxtaposition 
(e.g. for
objects one has $(+-)\otimes(-+-)=(+--+-)$ and for tangles $T_1\otimes 
T_2$ is
the tangle obtained by placing $T_2$ to the right of $T_1$), (B) a unit 
object
and unit morphism for this product, being the empty sequence and the empty
tangle respectively, and finally (C) an involution on objects given by
$+\leftrightarrow -$ in the sequences.

Thus the category of oriented tangles is very similar to the category of
1-cobordisms. There are however some important differences arising from 
the
fact that tangles are embedded in a 3-dimensional space, whereas the 
morphisms
of $1-{\bf Cobord}$ are abstract 1-dimensional manifolds. In the first
place, it is easy to find inequivalent tangles which are diffeomorphic as
abstract manifolds. Second, the product (A) is no longer commutative in 
${\bf OTa}$.

Tangles are very pleasing objects. They are a simultaneous generalization 
of both
braids and knots, since braids are tangles with all strands pointing
downwards, say, and knots, or more generally links, are just tangles from
$\emptyset$ to $\emptyset$. At the same time, tangles lend themselves 
well to
algebraization. Just as in the case of braids, where any braid can be 
written
as a composition of elementary braids, simple over- and undercrossings 
of two
adjacent strands, so also one can express any tangle as a composition 
of a
set  of elementary tangles. The equivalence of tangles under diffeomorphisms
of ${\bf R}^2\times [0,1]$ can be expressed algebraically as a set of
relations between these generators. Geometrically these correspond to a 
number
of tangle moves, analogous to the three Reidemeister moves for knot 
diagrams.
We refer to \cite{tur} for further details.

Now the main and beautiful result of the Turaev paper \cite{tur} is the
construction of a functor $F$ from the category of oriented tangles to 
${\bf
Vect}(k)$, where $k$ is now ${\bf Z}[q,q^{-1}]$, the ring of Laurent
polynomials in $q$ with integer coefficients. Let $V$ be a 
finite-dimensional
vector space over $k$ and $V^\ast$ its dual. Then the functor $F$ assigns
$\emptyset$ to $k$, $+$ to $V$, $-$ to $V^\ast$ and sequences of $+$ and $-$
are sent to the corresponding tensor product. Acting  on morphisms, $F$ 
sends a tangle
$O_1
\stackrel{T}{\longrightarrow} O_2$ to a linear map
$F(O_1) \stackrel{F(T)}{\longrightarrow}F(O_2)$. For instance the tangle
depicted in Figure ~\ref{fig:tangle}, corresponds under $F$ to a linear map
from $V^\ast\otimes V^\ast\otimes V\otimes V^\ast \otimes V$ to 
$V^\ast\otimes
V\otimes V^\ast$. This is how the term ``operator invariant'' arises, since
$F$ assigns a linear operator to each tangle, and ths assignment does not
depend on the particular representative of the equivalence class.

It would take us too far here to go into the details of the construction of
$F$. The idea is to find operators corresponding to the elementary tangles
referred to above, since by the functoriality of $F$ one has for the
composition of two tangles $F(T_1\circ T_2) = F(T_1)\circ F(T_2)$. In
particular the operators for the  elementary over- and undercrossings are
obtained from a certain class of quantum R-matrices, i.e. matrices 
satisfying
the quantum Yang-Baxter equation. This equation is the algebraic 
counterpart of
the third Reidemeister move for knot diagrams. For further details we refer,
once again, to \cite{tur}.

Now for the special case of tangles which are knots or links, the 
functor $F$
yields a linear map from $k$ to $k$, which may be identified with an 
element of
$k$. This element is essentially the Homfly polynomial of the knot or 
link, a
knot polynomial which generalises both the Alexander-Conway and Jones
polynomials. Thus Turaev's construction can be viewed as a generalisation 
to
tangles of the Homfly polynomial.

The point of discussing this example is hopefully clear by now: the 
functor $F$
which assigns operator invariants to tangles bears a very strong 
resemblance to
a 1-dimensional TQFT.  However, by virtue of the fact that the 
cobordisms are
now embedded manifolds the structure of the functor $F$ is far 
richer than that
of the simple 1-dimensional TQFT discussed in Sections 1 and 2. 
We are thus led
to seek a new class of TQFT-like theories which can encompass 
embedded
cobordisms, a point of view which will be explored now in the 
final section.

\section{\hspace{-4mm}.\hspace{2mm}TQFT FOR EMBEDDED MANIFOLDS}

\hspace*{0.8cm}
In this section we shall outline some ideas on how to modify the 
definition of
TQFT given in Section 2 in order to describe embedded manifolds. 
Our first
suggestion attempts to capture the features of the tangle example 
discussed in
the previous section. A ``tangle-type TQFT'' will be a functor from 
$(d,n)-{\bf Cobord}$
to ${\bf Vect}(k)$, where  $(d,n)-{\bf Cobord}$ is, loosely speaking, 
the category of
$d$-dimensional cobordisms embedded in $n$ dimensions. 


To define the objects of  $(d,n)-{\bf Cobord}$ we start by defining basic
objects. Consider embeddings of closed $(d-1)$-dimensional manifolds in
$I\times {\bf R}^{n-2}$ such that the image is contained in the interior of
$I\times {\bf R}^{n-2}$. Let these be connected in the following sense: 
if the
embedding is described by a map $f:\Sigma \rightarrow I\times 
{\bf R}^{n-2}$, there do
not exist open sets
$U_1, U_2
\subset I\times {\bf R}^{n-2}$, which are topologically open discs, 
such that
$U_1\cap U_2=\emptyset$,
$f(\Sigma)\subset U_1\cup U_2$ and $f(\Sigma) \cap U_1\neq \emptyset 
\neq f(\Sigma)\cap
U_2$. For each isotopy class of connected embeddings we choose a standard
representative. These are the basic objects. All other objects are 
obtained by
taking products of the basic objects, where the product is juxtaposition 
of the
embeddings in an analogous fashion to the case of tangles. This, to our 
mind, is
the natural way to generalize the objects of ${\bf OTa}$ to higher 
dimensions: the
standard embedding of a single point in $I\times {\bf R}$ identified with
$[1/2, 3/2]\times {\bf R}$ maps the point to $(1,0)$. Of course, it may be
interesting to study more general situations where the objects are 
unrestricted
embeddings, as in the analysis of braid statistics for particles in 
the plane
\cite{goldin}. The morphisms of $(d,n)-{\bf Cobord}$ are compact 
$d$-dimensional
manifolds embedded in ${\bf R}^{n-1}\times [0,1]$, whose boundary 
components
lie exclusively in the top and bottom hyperplanes ${\bf R}^{n-1}\times 
\{0\}$ and
${\bf R}^{n-1}\times \{1\}$, and both top and bottom boundaries belong to 
the
objects of $(d,n)-{\bf Cobord}$. The embedded d-dimensional manifolds are
identified up to orientation-preserving diffeomorphisms of 
${\bf R}^{n-1}\times
[0,1]$, which restrict to the identity on the top and bottom hyperplanes, as
in the case of tangles.

We have already seen an example of a tangle-type TQFT, namely the functor 
$F$
from $(1,3)-{\bf Cobord}$ to ${\bf Vect}(k)$ from the previous section.
Going up one dimension, another tangle-type TQFT is a functor from 
$(2,4)-{\bf
Cobord}$ to ${\bf Vect}(k)$. Now the objects of $(2,4)-{\bf
Cobord}$ are standard embeddings of closed 1-dimensional manifolds in 
$I\times
{\bf R}^2$, i.e. representatives of knot and link classes. The morphisms are
2-dimensional surfaces embedded in ${\bf R}^3\times [0,1]$ with objects as
their top and bottom boundaries, up to identification. A special case 
consists
of the 2-knots with empty top and bottom boundaries. Such objects have been
studied by a number of authors (see for instance \cite{carter}). 

A second context where an embedded TQFT structure appears is parallel 
transport
in vector bundles. Let $M$ be a smooth manifold. We can define a 
category ${\bf
Path}(M)$  as follows: the objects of  ${\bf Path}(M)$ are points of
$M$ and the set of morphisms from $m_1$ to $m_2$ is the set of smooth paths
from $m_1$ to $m_2$, i.e. maps $\gamma :[0,1]\rightarrow M$ such that
$\gamma(0)=m_1$, $\gamma(1)=m_2$, which are constant in $[0,\epsilon[$ and
$]\epsilon,1]$ for some $0<\epsilon<1/2$, identified up to a suitable
equivalence relation. The composition of two paths 
$m_1\stackrel{\gamma}{\longrightarrow}m_2$ and 
 $m_2\stackrel{\gamma'}{\longrightarrow}m_3$ is 
$m_1\stackrel{\gamma'\circ\gamma}{\longrightarrow}m_3$, the obvious path 
which
follows $\gamma$ and then $\gamma'$ at double speed. The equivalence
relation is such that this composition is well-defined and associative 
on path
equivalence classes. We refer to \cite{CP} for a description of an 
appropriate
equivalence relation, namely rank-1 homotopy. We remark that the extra
condition requiring the paths to be constant at their endpoints guarantees 
the
smoothness of the composition.

Now, given a vector bundle $E\stackrel{\pi}{\longrightarrow}M$ over $M$
with connection $\nabla$, we can define a ``parallel-transport-type TQFT" 
as a
functor $F$ from ${\bf Path}(M)$ to ${\bf Vect}(k)$, given by $m\mapsto
\pi^{-1}(m)$ (the fibre over $m$) for objects, and
$(m_1\stackrel{\gamma}{\longrightarrow}m_2)\mapsto
(\pi^{-1}(m_1)\stackrel{F(\gamma)}{\longrightarrow}\pi^{-1}(m_2))$ for 
morphisms,
where $F(\gamma)$ is the isomorphism between the fibres over the endpoints
induced by parallel transport along $\gamma$.

A natural modification of the previous example occurs when $\nabla$ is flat.
Then we can replace ${\bf Path}(M)$ by the homotopy groupoid, regarded as a
category whose morphisms are  homotopy classes of paths on $M$. 

Thus we see that a number of interesting examples have an embedded TQFT
structure. A question which naturally arises is: do these theories come 
about
from some classical theory and, if so,  what is the corresponding 
topological
action? For the lowest-dimensional tangle-type TQFT a candidate for 
the action
is presumably some version of the Kontsevich integral \cite{Kont} 
\cite{bn}: suppose the
tangle is a braid (i.e. none of the strands doubles back on itself). 
Any plane $x_3=c$ for
$c\in [0,1]$ constant intersects the braid in a fixed number of points, 
say $n$, and
thus varying the plane from $x_3=0$ to $x_3=1$ we get a film of $n$ points
moving in the plane. In other words, a braid defines a path in the
configuration space for $n$ identical non-coincident particles ${\bf
C}^n\setminus \Delta$, where we have identified ${\bf R}^2$ with ${\bf
C}$ and $\Delta$ is the diagonal subset of ${\bf C}^n$, i.e. the union of
hyperplanes $\left\{(z_1,\ldots ,z_n)|z_i=z_j\right \}$ for $i\neq j$. 
We can
introduce a flat connection on ${\bf C}^n\setminus \Delta$, namely the
Knizhnik-Zamolodchikov connection $A_{KZ}$, and the flatness implies 
that the
parallel transport along the path in ${\bf C}^n\setminus \Delta$ is
invariant under homotopy. In terms of braids this means that the parallel
transport is the same for two braids which are related by an
orientation-preserving diffeomorphism of
${\bf C}\times [0,1]$, which is the identity on the top and bottom planes
and which preserves each horizontal plane $x_3=c$. Now the Kontsevich 
integral
construction can be extended to tangles, \cite{Le-Mu}, \cite{Kass-Tur} 
giving a
functorial assignment from tangles to certain vector spaces generated by
so-called chord diagrams. Thus the topological action for tangles should be
something of the form ${\rm tr}\int f_\ast A_{KZ}$ where $f$ maps from 
$[0,1]$ to
some modified configuration space, appropriate to tangles rather than 
braids. 

In conclusion, the study of tangles from a TQFT angle suggests a 
promising method of
extending the TQFT approach to embedded manifolds. If the codimension 
of the embedded
manifolds is not too large, one can expect a considerably richer 
structure for embedded
TQFTs compared to pure TQFTs for the same dimension. From the 
physical point of view, it
is interesting to note that for this class of theories, we have 
written down a quantum
theory straight away, without having started from a classical theory, 
e.g. in terms of
some classical Lagrangian. Indeed it seems likely that the classical 
action will have a
complicated structure, if the tangle example is anything to go by. 
Thus the study of
TQFTs may suggest some new understanding of the nature of quantum 
theories, with no
dynamics to cloud the issue, at least in the first instance.

\section{\hspace{-4mm}.\hspace{2mm}ACKNOWLEDGEMENTS}

\hspace*{0.8cm}
I am grateful to the organisers of this Workshop for giving me the 
opportunity to
present this material, and to Cecile DeWitt-Morette, Jerry Goldin, 
Mauro Spera, Twareque
Ali and Tom\'{a}\v{s} S\'{y}kora for their encouraging reactions. 
I would also like to
thank my fellow participants in the Lisbon ``TQFT Club'' meetings for 
many exchanges and
moral support. This work was supported by Junta Nacional de 
Investiga\c{c}\~ao
Cient\'{\i}fica e Tecnol\'{o}gica and the Praxis XXI programme, as well 
as by the European
Union under the Human Capital and Mobility programme.


\begin{thebibliography}{99}

\footnotesize

\bibitem{at} M. Atiyah, {\em Publ. Math. Inst. Hautes \'{E}tudes Sci. 
(Paris)} 68, 175,
(1989)

\bibitem{yetter} D. N. Yetter, Triangulations and TQFT's, {\em in} 
``Quantum Topology'',
L. H. Kauffman and R. A. Baadhio eds., World Scientific, Singapore (1993)

\bibitem{baezdolan} J. Baez and J. Dolan, {\em J.Math. Phys.} 36, 6073, 
(1995)

\bibitem{crane} L. Crane, {\em J.Math. Phys.} 36, 6180, (1995)

\bibitem{sawin1} S.Sawin, {\em J.Math. Phys.} 36, 6673, (1995)

\bibitem{qft+jones} E. Witten, {\em Commun. Math. Phys.} 121, 351, (1989)

\bibitem{turbook} V. G. Turaev, ``Quantum Invariants of Knots and 
3-Manifolds'', de
Gruyter, Berlin (1994) 

\bibitem{blauthomp} M. Blau and G. Thompson, {\em J.Math. Phys.} 36, 
2192, (1995)

\bibitem{masbaum} C. Blanchet, N. Habegger, G. Masbaum and P. Vogel, 
{\em Topology} 34,
883 (1995)

\bibitem{sawin2}  S.Sawin, {\em Bull. Amer. Math. Soc.} 33, 413, (1996)

\bibitem{freed} D. S. Freed, Lectures on topological quantum field 
theory, {\em in}
``Integrable Systems, Quantum Groups, and Quantum Field Theory'', 
L. A. Ibort and M. A.
Rodriguez eds., Kluwer, Amsterdam (1993)

\bibitem{tur} V. G. Turaev, {\em Math. USSR Izvestija} 35, 411, (1990)

\bibitem{goldin} G. A. Goldin and D. H. Sharp, {\em Phys. Rev. Lett.} 
76, 1183, (1996)

\bibitem{carter} J. Carter and M. Saito, Knotted surfaces, braid movies, 
and
beyond, {\em in} ``Knots and Quantum Gravity'', J. Baez ed., 
Oxford University Press,
Oxford (1994) 

\bibitem{CP} A. Caetano and R. F. Picken, {\em Int. J. Math.} 5, 835, (1994) 

\bibitem{Kont} M. Kontsevich, {\em Adv. Sov. Math.} 16, 137, (1993)

\bibitem{bn} D. Bar-Natan,{\em Topology} 34, 423 (1995)

\bibitem{Le-Mu} T. Q. T. Le and J. Murakami, ``Representation of the 
category of tangles
by Kontsevich's iterated integral'' {\em Commun. Math. Phys.} (to appear)

\bibitem{Kass-Tur} C. Kassel and V. Turaev, ``Chord diagram invariants 
of tangles and
graphs'', {\em Preprint} (1996)







\end{thebibliography}
\end{document}